\documentclass[11pt,a4paper,reqno]{amsart}
\usepackage{amsmath,amssymb}
\usepackage[a4paper,left=79pt,right=79pt,top=80pt,bottom=80pt]{geometry}
\usepackage[hidelinks]{hyperref}
\hypersetup{
  pdftitle={Uniqueness Sequences in the Bargmann--Fock Space},
  pdfauthor={Shengzhao Hou and Xiangrun Yang},
  pdfsubject={Jensen potentials and uniqueness sequences in Fock space}
}
\makeatletter
\renewcommand{\subsection}{%
  \@startsection{subsection}{2}{\z@}%
    {.5\linespacing\@plus.7\linespacing}{.35\linespacing}%
    {\normalfont\bfseries\raggedright}}
\makeatother
\allowdisplaybreaks[1]
\newcommand{\C}{\mathbb C}
\newcommand{\R}{\mathbb R}
\newcommand{\N}{\mathbb N}
\newcommand{\Z}{\mathbb Z}
\newcommand{\Q}{\mathbb Q}
\newcommand{\T}{\mathbb T}
\newcommand{\F}{\mathcal F^2}
\newcommand{\J}{\mathcal J}
\newcommand{\E}{\mathcal E}
\newcommand{\Hol}{\operatorname{Hol}}
\newcommand{\SH}{\operatorname{SH}}
\newcommand{\supp}{\operatorname{supp}}
\newcommand{\dist}{\operatorname{dist}}
\newcommand{\ord}{\operatorname{ord}}
\newcommand{\im}{\operatorname{Im}}
\newcommand{\Ree}{\operatorname{Re}}
\newcommand{\ind}{\mathbf 1}
\newcommand{\dd}{\,d}
\newcommand{\norm}[1]{\left\|#1\right\|}
\newtheorem{theorem}{Theorem}[section]
\newtheorem{lemma}[theorem]{Lemma}
\newtheorem{proposition}[theorem]{Proposition}
\newtheorem{corollary}[theorem]{Corollary}
\theoremstyle{remark}
\newtheorem{remark}[theorem]{Remark}
\numberwithin{equation}{section}
\title[Uniqueness Sequences in the Bargmann--Fock Space]{Uniqueness Sequences in the Bargmann--Fock Space}
\author{Shengzhao Hou}
\author{Xiangrun Yang}
\thanks{School of Mathematical Sciences, Soochow University, Suzhou 215006, China.}
\thanks{Corresponding author: Shengzhao Hou, \href{mailto:shou@suda.edu.cn}{shou@suda.edu.cn}.}
\thanks{This work was supported by the National Natural Science Foundation of China (Grant No. 12371133).}
\date{}
\begin{document}
\begin{abstract}
We characterize uniqueness sequences in the Bargmann--Fock space by the divergence of a functional involving Jensen potentials and Gaussian relative entropy. The criterion applies to every locally finite sequence with prescribed multiplicities. Its proof combines an exact subharmonic duality with holomorphic realization at the same Gaussian weight. The radial restriction yields an integral criterion that is necessary and sufficient for unions of rotation orbits satisfying a bounded-overlap condition. We also show that every pairwise distinct angular rearrangement of the nonzero radii of a lattice through the origin with fundamental area $\pi$ is a uniqueness sequence. In particular, independently rotated square-lattice radii give a complete coherent-state system almost surely at the critical scale. Finally, we reduce the general criterion to countably many polynomial tests and give finite approximations with explicit error bounds.
\end{abstract}
\maketitle
\enlargethispage{2pt}

\section*{Introduction}
Let $A$ be planar Lebesgue measure and let $\dd\gamma(z)=\pi^{-1}e^{-|z|^2}\dd A(z)$. We equip the Bargmann--Fock space
\[
\F=\left\{f\in\Hol(\C):\norm{f}_{\F}^{2}
=\int_\C |f(z)|^2\dd\gamma(z)<\infty\right\}
\]
with its usual Hilbert space norm. A locally finite sequence $\Lambda$ is a uniqueness sequence for $\F$ if every function in $\F$ that vanishes on $\Lambda$, with the prescribed multiplicities, is identically zero. Additional zeros are allowed. An exact zero sequence, by definition, is the complete zero sequence of a nonzero function in $\F$. The two notions arise naturally in the study of zeros and reproducing kernels in Fock spaces (see \cite{AadiBouyaOmari,AadiOmari,Zhu}).

In joint work with P.~Ma, the first author characterized exact zero sequences by a relative-entropy formula for a renormalized canonical-product potential, subject to a quadratic counting condition and two complex moment constraints \cite[Theorem~2.1]{HouMa}. In the present uniqueness problem, additional zeros are allowed. We therefore minimize over subharmonic functions whose Riesz measures dominate the prescribed counting measure. Our criterion applies to every locally finite sequence. Jensen measures and subharmonic minorants have also been used in zero-subset problems for weighted spaces of holomorphic functions \cite{Khabibullin}.

Theorem~\ref{thm:unique} characterizes uniqueness by the divergence of a supremum over smooth Jensen densities. The score of a density is the sum of twice its Jensen potential over the prescribed zeros, counted with multiplicity, minus its relative entropy with respect to $\gamma$. Proposition~\ref{prop:duality} identifies the exponential of this supremum with the minimum of a normalized subharmonic integral when the value is finite. Proposition~\ref{prop:realize} supplies a nonzero entire function satisfying the prescribed vanishing conditions with finite norm at the same Gaussian weight.

Section~\ref{sec:radial} evaluates the supremum over radial densities. This formula recovers the Jensen integral obstruction to a Fock function vanishing on a prescribed sequence. Divergence is also necessary for uniqueness for unions of regular rotation orbits under a bounded-overlap condition. For lattices through the origin with fundamental area $\pi$, a uniform estimate for the integrated radial count shows that every pairwise distinct angular rearrangement of the nonzero lattice radii is a uniqueness sequence. This resolves the critical case of the independent square-radius model posed by Kononova \cite[Introduction, item~2(c)]{KononovaRandom}. For general critical configurations, even the full radial counting function together with finitely many angular moments can fail to determine uniqueness \cite[Theorem~8.7]{HouMa}.

Section~\ref{sec:polynomials} reduces the full supremum to tests indexed by rational polynomials and rational smoothing radii. For each test we give a finite approximation by disk arrangements and bound its error explicitly. If the sequence can be listed with certified coordinates in every test disk, the resulting lower bounds increase to the variational supremum.

\section{A variational criterion for uniqueness sequences}\label{sec:general}
\subsection{Notation and statement of the theorem}
Let $\Lambda=(\lambda_j)_{j\in I}$ be a locally finite sequence in $\C$, where $I$ is finite or countable. Local finiteness means that every compact set contains only finitely many terms, including repetitions. Write $\operatorname{ran}\Lambda$ for the set of distinct points in the sequence. For $\lambda\in\operatorname{ran}\Lambda$, let $m_\lambda=\#\{j\in I:\lambda_j=\lambda\}$, and set
\[
n_\Lambda=\sum_{\lambda\in\operatorname{ran}\Lambda}m_\lambda\delta_\lambda,
\qquad
N_\Lambda=\{f\in\F:f^{(k)}(\lambda)=0
\text{ for }\lambda\in\operatorname{ran}\Lambda,\ 0\leq k<m_\lambda\}.
\]
All sums weighted by $m_\lambda$ run over distinct points. If $f$ is a nonzero entire function, write
\[
n_f=\sum_{f(a)=0}\ord_a(f)\,\delta_a.
\]
For $f\in\F$ with $f\not\equiv0$, membership in $N_\Lambda$ is equivalent to $n_f\geq n_\Lambda$. Thus uniqueness means $N_\Lambda=\{0\}$, while exact realization means $n_f=n_\Lambda$ for some nonzero $f\in\F$.

For a subharmonic function $u$, let
\[
\nu_u=\frac1{2\pi}\Delta u.
\]
In particular, $\nu_{\log|f|}=n_f$. Subharmonic functions are always taken with their upper semicontinuous representatives and are not identically $-\infty$.

We first assume that the sequence avoids the origin, that is, $\lambda_j\ne0$ for every $j\in I$. For a compactly supported probability density $\rho$, define
\begin{equation}\label{eq:potential}
V_\rho(w)=\int_\C\log|w-z|\rho(z)\dd A(z)-\log|w|,
\qquad w\ne0,
\end{equation}
whenever the integral is defined. Our fixed family of smooth tests is
\begin{equation}\label{eq:J}
\J=\left\{\rho\in C_c^\infty(\C):\rho\geq0,\quad
\int_\C\rho\dd A=1,\quad V_\rho(w)\geq0\text{ for }w\ne0\right\}.
\end{equation}
These are smooth Jensen densities at the origin. Their Gaussian relative entropy is
\begin{equation}\label{eq:entropy}
H_\gamma(\rho)=\int_\C\rho(z)\log\bigl(\pi\rho(z)e^{|z|^2}\bigr)\dd A(z),
\qquad0\log0=0.
\end{equation}
We use the same formula for other probability densities, allowing the value $+\infty$. Define
\begin{equation}\label{eq:E}
\E(\Lambda)=\sup_{\rho\in\J}
\left\{2\sum_{\lambda}m_\lambda V_\rho(\lambda)-H_\gamma(\rho)\right\}.
\end{equation}
Each smooth test has finite entropy. By Lemma~\ref{lem:jensen}, its potential has bounded support, so the sum is finite.

\begin{theorem}[General criterion]\label{thm:unique}
Let $\Lambda=(\lambda_j)_{j\in I}$ be a locally finite sequence in $\C\setminus\{0\}$, with repetitions allowed. Then
\[
N_\Lambda=\{0\}\quad\Longleftrightarrow\quad\E(\Lambda)=+\infty.
\]
In particular, $\Lambda$ is a uniqueness sequence precisely when the Jensen--entropy functional diverges.
\end{theorem}

To prove Theorem~\ref{thm:unique}, we first establish the subharmonic duality in Proposition~\ref{prop:duality}. Proposition~\ref{prop:realize} then gives the required holomorphic function at the original Gaussian weight, using the weighted Fock dimension theorem of Borichev, Le, and Youssfi \cite{BorichevLeYoussfi}.

For a sequence containing the origin, choose $a\notin\operatorname{ran}\Lambda$. The operator
\begin{equation}\label{eq:mag}
(T_a f)(z)=e^{-\overline a z-|a|^2/2}f(z+a)
\end{equation}
is a surjective isometry of $\F$ carrying $N_\Lambda$ onto $N_{\Lambda-a}$. Theorem~\ref{thm:unique} therefore applies to $\Lambda-a$. The equivalence with uniqueness shows that finiteness of the supremum is independent of the choice of $a$.

\subsection{Jensen potentials and entropy}
We use the standard facts about subharmonic functions and Riesz measures collected in \cite{Ransford}.
\begin{lemma}\label{lem:jensen}
Let $\rho\in\J$, and suppose $\supp\rho\subset D(0,R)$. Then $V_\rho=0$ on $\{|w|>R\}$. If $u$ is subharmonic on a neighborhood of $\overline{D(0,S)}$, where $S>R$, and $u(0)>-\infty$, then
\begin{equation}\label{eq:jensen}
\int_\C u\rho\dd A-u(0)=\int_{D(0,S)}V_\rho\dd\nu_u.
\end{equation}
In particular, $\int u\rho\dd A\geq u(0)$.
\end{lemma}
\begin{proof}
The circular mean formula is
\begin{equation}\label{eq:circular-mean}
\frac1{2\pi}\int_0^{2\pi}\log|re^{it}-z|\dd t
=\log\max\{r,|z|\},\qquad z\in\C,\quad r>0.
\end{equation}
For $|z|<r$, this follows by taking the real part of the mean value identity for
$\log(1-z e^{-it}/r)$; the case $|z|>r$ follows by interchanging the roles of $r$ and $|z|$, and the boundary case follows by continuity.
For $r>R$, Fubini's theorem and \eqref{eq:circular-mean} give
\[
\frac1{2\pi}\int_0^{2\pi}V_\rho(re^{it})\dd t=0.
\]
Since $V_\rho$ is continuous and nonnegative away from the origin, it vanishes for $|w|>R$.

Choose $R<r_1<r_2<S$. The local Riesz decomposition gives
\[
u(z)=h(z)+\int_{D(0,r_2)}\log|z-w|\dd\nu_u(w),\qquad z\in D(0,r_1),
\]
with $h$ harmonic. We regard $V_\rho$ as a locally integrable function across the origin. Since
\[
\Delta V_\rho=2\pi(\rho\dd A-\delta_0),\qquad
\supp\Delta V_\rho\subset\supp V_\rho,
\]
and $\rho\dd A$ has no atom at zero, we have
$\supp\rho\cup\{0\}\subset\supp V_\rho$.
Choose $\chi\in C_c^\infty(D(0,r_1))$ equal to one near $\supp V_\rho$.
Distributional integration by parts yields
\[
\int h\rho\dd A-h(0)
=\frac1{2\pi}\langle V_\rho,\Delta(\chi h)\rangle=0.
\]
Here $\chi=1$ on $\supp\rho\cup\{0\}$, and $\Delta(\chi h)$ vanishes near $\supp V_\rho$.
Applying Fubini's theorem to the logarithmic part of the decomposition now gives
\[
\int u\rho\dd A-u(0)
=\int_{D(0,r_2)}V_\rho\dd\nu_u
=\int_{D(0,S)}V_\rho\dd\nu_u.
\]
The interchange is justified by $u(0)>-\infty$, which controls the negative logarithmic singularity at the origin, and by
\[
\sup_{w\in D(0,r_2)}\int_\C|\log|z-w||\rho(z)\dd A(z)<\infty.
\]
This proves the identity and its nonnegativity.
\end{proof}

In particular, $\rho\dd A$ is a Jensen measure for subharmonic functions on any centered disk containing its support compactly. Apply the lemma on an intermediate disk when $u(0)>-\infty$; the inequality is automatic when $u(0)=-\infty$. Applying it to both signs of a harmonic function gives equality. In particular,
\begin{equation}\label{eq:barycenter}
\int_\C z\rho(z)\dd A(z)=0.
\end{equation}

\begin{lemma}\label{lem:entropy}
Suppose that $u\in\SH(\C)$ and $u(0)>-\infty$. For every $\rho\in\J$,
\begin{equation}\label{eq:entropyineq}
2\int u\rho\dd A-H_\gamma(\rho)\leq\log\int e^{2u}\dd\gamma.
\end{equation}
Consequently, if $\nu_u\geq n_\Lambda$, then
\begin{equation}\label{eq:entropyE}
\E(\Lambda)\leq\log\int e^{2u}\dd\gamma-2u(0).
\end{equation}
\end{lemma}
\begin{proof}
The assertion is immediate if $Z:=\int e^{2u}\dd\gamma=\infty$.
Otherwise put $q=Z^{-1}e^{2u}\pi^{-1}e^{-|z|^2}$ and $G=\{\rho>0\}$.
Jensen's inequality for the concave logarithm gives
\[
\int_G\rho\log\frac q\rho\dd A
\leq\log\int_G q\dd A\leq0.
\]
Expanding the logarithm yields
\[
0\leq H_\gamma(\rho)-2\int u\rho\dd A+\log Z.
\]
All logarithmic integrals on the left are finite by local integrability of $u$ and the smooth compact support of $\rho$. This is the usual relative entropy inequality; see \cite{CoverThomas}.
If $\nu_u\geq n_\Lambda$, Lemma~\ref{lem:jensen} gives
\[
\sum_{\lambda\in\operatorname{ran}\Lambda}m_\lambda V_\rho(\lambda)
\leq\int V_\rho\dd\nu_u=\int u\rho\dd A-u(0).
\]
Combining the two inequalities and taking the supremum proves \eqref{eq:entropyE}.
\end{proof}
Every smooth, compactly supported, radial probability density $\rho(z)=\varphi(|z|)$ belongs to $\J$, since the circular mean identity gives
\[
V_\rho(w)=2\pi\int_{|w|}^{\infty}\varphi(r)r\log\frac r{|w|}\dd r\geq0,\qquad w\neq0.
\]
Thus $\J$ is nonempty.

To show that $\E(\Lambda)\geq0$, choose smooth radial cutoffs $0\leq\chi_R\leq1$ equal to one on $D(0,R)$, with compact support, and put
\[
c_R=\int_\C\chi_R\dd\gamma,\qquad
\rho_R(z)=c_R^{-1}\chi_R(z)\pi^{-1}e^{-|z|^2}.
\]
Then $\rho_R\in\J$, $c_R\to1$, and
\[
0\leq H_\gamma(\rho_R)
=-\log c_R+c_R^{-1}\int_\C\chi_R\log\chi_R\dd\gamma
\leq-\log c_R\longrightarrow0.
\]
Since $V_{\rho_R}\geq0$, we obtain $\E(\Lambda)\geq-H_\gamma(\rho_R)\to0$. Moreover, $H_\gamma(\rho)\geq0$ for all $\rho\in\J$, so $\E(\varnothing)=0$.

\subsection{The subharmonic variational problem}
Write $D_R=D(0,R)$ and let $\Lambda_R$ be the subsequence consisting of all terms $\lambda_j$ with $|\lambda_j|<R$, retaining repetitions. Define
\[
P_R(z)=\prod_{\lambda\in\operatorname{ran}\Lambda_R}(1-z/\lambda)^{m_\lambda},\qquad q_R=\log|P_R|.
\]
The empty product is one. Thus $q_R(0)=0$ and $\nu_{q_R}=n_{\Lambda_R}$. Consider the local minimum
\begin{equation}\label{eq:ZR}
Z_R=\inf\left\{\int_{D_R}e^{2u}\dd\gamma:u\in\SH(D_R),\ u(0)\geq0,\ \nu_u\geq n_{\Lambda_R}\right\}.
\end{equation}

\begin{lemma}\label{lem:local}
For every $R>0$, we have $0<Z_R<\infty$, and the infimum is attained by a function $u_R$ satisfying $u_R(0)=0$.
\end{lemma}
\begin{proof}
The function $q_R$ is admissible and has a finite exponential integral.
For $0<s<R$, the submean inequality and Jensen's inequality give, for every admissible $u$,
\[
\int_{D_R}e^{2u}\dd\gamma
\geq e^{-s^2}s^2
\exp\left(\frac1{\pi s^2}\int_{D_s}2u\dd A\right)
\geq e^{-s^2}s^2>0.
\]
Thus $0<Z_R<\infty$.

Let $(u_j)$ be a minimizing sequence. For each compact $K\subset D_R$, choose $s>0$ with $D(z,s)\subset D_R$ for $z\in K$.
Since $e^{2u_j}$ is subharmonic,
\[
e^{2u_j(z)}\leq\frac{e^{R^2}}{s^2}\int_{D_R}e^{2u_j}\dd\gamma,\qquad z\in K.
\]
Hence $(u_j)$ is locally bounded above. The conditions $u_j(0)\geq0$ exclude convergence to $-\infty$.
Subharmonic compactness therefore gives a subsequence converging in $L^1_{\mathrm{loc}}(D_R)$ to a subharmonic function $u_R$.
Distributional convergence preserves $\nu_{u_R}\geq n_{\Lambda_R}$.
For each $s<R$,
\[
\frac1{\pi s^2}\int_{D_s}u_R\dd A
=\lim_{j\to\infty}\frac1{\pi s^2}\int_{D_s}u_j\dd A\geq0.
\]
Letting $s\downarrow0$ yields $u_R(0)\geq0$.
After taking a further almost everywhere convergent subsequence, Fatou's lemma proves attainment.
If $u_R(0)>0$, subtracting $u_R(0)$ strictly decreases the integral while preserving admissibility. Thus $u_R(0)=0$.
\end{proof}

\begin{lemma}\label{lem:min-density}
Let $u_R$ be the minimizer given by Lemma~\ref{lem:local}, and define the probability density
\begin{equation}\label{eq:pR}
p_R(z)=Z_R^{-1}e^{2u_R(z)}\pi^{-1}e^{-|z|^2}\ind_{D_R}(z).
\end{equation}
Then $V_{p_R}(w)\geq0$ for $w\neq0$. Moreover,
\begin{equation}\label{eq:pR-entropy}
\begin{aligned}
H_\gamma(p_R)=2\int_\C u_Rp_R\dd A-\log Z_R&<\infty,\\
2\int_\C q_Rp_R\dd A-H_\gamma(p_R)&\geq\log Z_R.
\end{aligned}
\end{equation}
\end{lemma}
\begin{proof}
Let $v$ be smooth and subharmonic on $\C$, with $v(0)=0$.
For $t\geq0$, the function $u_R+tv$ is admissible.
Since $v$ is bounded on $D_R$, differentiation at $t=0^+$ gives
\begin{equation}\label{eq:variation}
\int_\C vp_R\dd A\geq0.
\end{equation}
Apply this to
\[
v_{w,\varepsilon}(z)=\tfrac12\log(|z-w|^2+\varepsilon^2)
-\tfrac12\log(|w|^2+\varepsilon^2),\qquad w\ne0.
\]
As $\varepsilon\downarrow0$, monotone convergence, after subtraction from a common upper bound on $D_R$, gives
\[
\int_\C\log|z-w|p_R(z)\dd A(z)\geq\log|w|.
\]
In particular this integral is not $-\infty$, and $V_{p_R}(w)\geq0$.

The logarithmic singularities of $q_R$ lie in $D_R$, where $p_R$ is locally bounded, so $q_R\in L^1(p_R\dd A)$.
Set $a=q_R-u_R$, with differences understood almost everywhere. The admissible convex combinations
$u_t=(1-t)u_R+tq_R$ have finite exponential integrals by H\"older's inequality, and minimality gives
\begin{equation}\label{eq:expvariation}
\int_\C e^{2ta}p_R\dd A\geq1,\qquad 0<t\leq1.
\end{equation}
Choose $M$ such that $q_R\leq M$ on $D_R$.
On $\{a>0\}$ and for $0<t\leq\tfrac12$,
\[
0\leq\frac{e^{2ta}-1}{t}\leq2ae^{2ta},
\qquad
 ae^{2ta}e^{2u_R}
\leq e^{2M}(M-u_R)e^{-(M-u_R)}\leq e^{2M-1}.
\]
Consequently the positive part of the difference quotient is dominated, after multiplication by $p_R$, by the integrable constant
$2e^{2M-1}/(\pi Z_R)$ on $D_R$.
It converges in integral to $2\int a_+p_R\dd A<\infty$.
On $\{a<0\}$, the nonnegative functions
$(1-e^{2ta})/t$ increase to $2a_-$ as $t\downarrow0$.
By \eqref{eq:expvariation}, their integrals are bounded by those of the positive difference quotients.
Monotone convergence therefore gives
\[
a\in L^1(p_R\dd A),\qquad \int_\C a p_R\dd A\geq0.
\]
Since $q_R$ is integrable, so is $u_R$.
Expanding the logarithm in the entropy then gives
\[
H_\gamma(p_R)=2\int_\C u_Rp_R\dd A-\log Z_R<\infty,
\]
and
\[
2\int_\C q_Rp_R\dd A-H_\gamma(p_R)
=\log Z_R+2\int_\C a p_R\dd A\geq\log Z_R.
\]
\end{proof}

The next lemma allows us to use the minimizing density $p_R$ even without a bound near $\partial D_R$.

\begin{lemma}\label{lem:smoothing}
Let $p$ be a compactly supported probability density with $V_p(w)\geq0$ for $w\neq0$ and $H_\gamma(p)<\infty$. Let $\eta_\varepsilon$ be a smooth radial probability density supported in $D(0,\varepsilon)$. Set
\[
p_\varepsilon(z)=(p*\eta_\varepsilon)(z)
=\int_\C p(u)\eta_\varepsilon(z-u)\dd A(u).
\]
Then $p_\varepsilon\in\J$, and
\begin{equation}\label{eq:smoothing}
V_{p_\varepsilon}(w)\geq V_p(w),\qquad
H_\gamma(p_\varepsilon)\leq H_\gamma(p)+s_\varepsilon,\qquad
s_\varepsilon=\int|z|^2\eta_\varepsilon(z)\dd A(z)\leq\varepsilon^2.
\end{equation}
\end{lemma}
\begin{proof}
For every $a\in\C$, the circular mean formula gives
\[
\int_\C\log|a-v|\eta_\varepsilon(v)\dd A(v)\geq\log|a|.
\]
Integrating this inequality against $p(u)\dd A(u)$ with $a=w-u$ yields
$V_{p_\varepsilon}(w)\geq V_p(w)$.
The positive parts of the logarithms are bounded on the relevant compact sets, while integration against the smooth kernel controls their negative parts, which justifies Fubini's theorem.
Thus $p_\varepsilon$ is a smooth, compactly supported Jensen density.

Convexity of $x\log x$ and translation invariance of area give
\begin{equation}\label{eq:smoothing-convexity}
\int_\C p_\varepsilon\log p_\varepsilon\dd A
\leq\int_\C p\log p\dd A.
\end{equation}
Since $\eta_\varepsilon$ is radial, it has mean zero. Hence
\[
\int_\C|z|^2p_\varepsilon(z)\dd A(z)
=\int_\C|u|^2p(u)\dd A(u)+s_\varepsilon.
\]
Combining the last two statements proves the entropy estimate.
The support condition on $\eta_\varepsilon$ gives $s_\varepsilon\leq\varepsilon^2$.
\end{proof}

\begin{lemma}\label{lem:ZR}
For every $R>0$,
\begin{equation}\label{eq:ZRbound}
Z_R\leq e^{\E(\Lambda)}.
\end{equation}
\end{lemma}
\begin{proof}
Apply Lemma~\ref{lem:smoothing} to $p_R$. Nonnegativity of the potentials, followed by Lemma~\ref{lem:min-density}, gives
\begin{align*}
\E(\Lambda)&\geq 2\sum_{\lambda\in\operatorname{ran}\Lambda}m_\lambda V_{p_R*\eta_\varepsilon}(\lambda)
-H_\gamma(p_R*\eta_\varepsilon)\\
&\geq2\sum_{\lambda\in\operatorname{ran}\Lambda_R}m_\lambda V_{p_R}(\lambda)
-H_\gamma(p_R)-s_\varepsilon\\
&=2\int_\C\biggl(\sum_{\lambda\in\operatorname{ran}\Lambda_R}m_\lambda
 \log\left|1-\frac z\lambda\right|\biggr)p_R(z)\dd A(z)
-H_\gamma(p_R)-s_\varepsilon\\
&=2\int_\C q_Rp_R\dd A-H_\gamma(p_R)-s_\varepsilon\\
&\geq\log Z_R-s_\varepsilon.
\end{align*}
Letting $\varepsilon\downarrow0$, we obtain the desired inequality.
\end{proof}

\begin{proposition}[Jensen--entropy duality]\label{prop:duality}
Let $\Lambda=(\lambda_j)_{j\in I}$ be a locally finite sequence in $\C\setminus\{0\}$, where $I$ is finite or countable and repetitions are allowed. For each term value $\lambda$, let $m_\lambda$ be its number of occurrences. Define
\begin{equation}\label{eq:U}
\begin{gathered}
U_\Lambda=\{u\in\SH(\C):u(0)=0,\ \nu_u\geq n_\Lambda\},\\
R(\Lambda)=\inf_{u\in U_\Lambda}\int_\C e^{2u}\dd\gamma.
\end{gathered}
\end{equation}
Then
\begin{equation}\label{eq:duality}
R(\Lambda)=e^{\E(\Lambda)},
\end{equation}
where both sides may be infinite. If $R(\Lambda)<\infty$, the infimum in \eqref{eq:U} is attained by a unique subharmonic function.
\end{proposition}
\begin{proof}
The class $U_\Lambda$ is nonempty: by the Weierstrass factorization theorem, there is an entire function $F$ with zero divisor $n_\Lambda$. Since $0\notin\Lambda$, we may normalize $F(0)=1$, and then $\log|F|\in U_\Lambda$. This does not imply that $R(\Lambda)$ is finite. By \eqref{eq:entropyE},
\begin{equation}\label{eq:Ebound}
\E(\Lambda)\leq\log R(\Lambda).
\end{equation}
Suppose first that $\E(\Lambda)<\infty$, and choose $R_j\uparrow\infty$. The minimizers satisfy $u_{R_j}(0)=0$ and
\[\int_{D_{R_j}}e^{2u_{R_j}}\dd\gamma\leq e^{\E(\Lambda)}.\]
For $z\in D_s$ and $R_j>s+1$, the submean inequality on $D(z,1)$ gives
\[
e^{2u_{R_j}(z)}\leq e^{(s+1)^2}\int_{D_{R_j}}e^{2u_{R_j}}\dd\gamma.
\]
Hence the sequence is locally bounded above. Since $u_{R_j}(0)=0$, subharmonic compactness and a diagonal argument yield a subsequence converging in $L^1_{\mathrm{loc}}(\C)$ and almost everywhere to $u\in\SH(\C)$. By the submean inequality and $u_{R_j}(0)=0$, the averages of $u$ over disks centered at the origin are nonnegative. Letting the radii tend to zero gives $u(0)\geq0$. Distributional convergence preserves $\nu_u\geq n_\Lambda$, since every compact set is eventually contained in $D_{R_j}$. For fixed $s>0$ and all sufficiently large $j$, Fatou's lemma gives
\[
\int_{D_s}e^{2u}\dd\gamma
\leq\liminf_{j\to\infty}\int_{D_s}e^{2u_{R_j}}\dd\gamma
\leq e^{\E(\Lambda)}.
\]
Letting $s\uparrow\infty$ yields $\int_\C e^{2u}\dd\gamma\leq e^{\E(\Lambda)}$.
Replacing $u$ by $u-u(0)$ provides a member of $U_\Lambda$ with no larger integral. Combining this with \eqref{eq:Ebound} proves equality and attainment. If $\E(\Lambda)=+\infty$, \eqref{eq:entropyE} rules out a finite integral for any member of $U_\Lambda$, so \eqref{eq:duality} holds in that case as well.

Finally, $U_\Lambda$ is convex. If two distinct minimizers have finite integral, strict convexity of $e^{2x}$ gives a strictly smaller integral at their midpoint, unless they agree almost everywhere. Subharmonic representatives agreeing almost everywhere agree everywhere. This proves uniqueness.
\end{proof}

\begin{corollary}\label{cor:exhaust}
The quantities $Z_R$ are nondecreasing in $R$, and
\[\lim_{R\to\infty}Z_R=R(\Lambda)=e^{\E(\Lambda)}.\]
\end{corollary}
\begin{proof}
Restriction to smaller disks shows that $Z_R$ is nondecreasing in $R$. Moreover, restricting any globally admissible function to $D_R$ gives $Z_R\leq R(\Lambda)$. Hence
\[
L:=\lim_{R\to\infty}Z_R\leq R(\Lambda).
\]
Suppose $L<\infty$, and choose $R_j\uparrow\infty$. By the preceding compactness argument, a subsequence of $(u_{R_j})$ converges in $L^1_{\mathrm{loc}}(\C)$ and almost everywhere to a subharmonic function $u$ satisfying $u(0)\geq0$ and $\nu_u\geq n_\Lambda$. For fixed $s$, Fatou's lemma on $D_s$ gives $\int_{D_s}e^{2u}\dd\gamma\leq L$. Letting $s\uparrow\infty$, we obtain
\[
\int_\C e^{2u}\dd\gamma\leq L.
\]
Since $v=u-u(0)$ is globally admissible, we obtain
\[
R(\Lambda)\leq\int_\C e^{2v}\dd\gamma
\leq\int_\C e^{2u}\dd\gamma\leq L.
\]
If $L=\infty$, the inequality $Z_R\leq R(\Lambda)$ already implies $R(\Lambda)=\infty$. Thus
\[
R(\Lambda)=\lim_{R\to\infty}Z_R.
\]
\end{proof}

\subsection{Holomorphic realization at the Gaussian weight}
For $r>0$, let
\[
M_ru(z)=\frac1{\pi r^2}\int_{D(z,r)}u(w)\dd A(w).
\]
The following estimate supplies the weighted control needed for the realization argument. For a circular-average estimate for logarithms of entire functions, see \cite[Lemma~2.2]{HouMa}.

\begin{lemma}\label{lem:gaussavg}
For every real $u\in L^1_{\mathrm{loc}}(\C)$ and $r>0$,
\begin{equation}\label{eq:gaussavg}
\int_\C e^{2M_ru}\dd\gamma\leq e^{r^2}\int_\C e^{2u}\dd\gamma.
\end{equation}
\end{lemma}
\begin{proof}
The function $w\mapsto\Ree(\overline z(w-z))$ has mean zero on $D(z,r)$. Jensen's inequality therefore gives
\begin{align*}
e^{2M_ru(z)-|z|^2}
&\leq\frac1{\pi r^2}\int_{D(z,r)}e^{2u(w)-2\Ree(\overline z(w-z))-|z|^2}\dd A(w)\\
&=\frac1{\pi r^2}\int_{D(z,r)}e^{2u(w)-|w|^2+|w-z|^2}\dd A(w)\\
&\leq\frac{e^{r^2}}{\pi r^2}\int_{D(z,r)}e^{2u(w)-|w|^2}\dd A(w).
\end{align*}
Integrating in $z$ and using Tonelli's theorem proves \eqref{eq:gaussavg}.
\end{proof}

For a positive locally finite measure $\nu$, define its integer atomic part and its remainder by
\begin{equation}\label{eq:atomic}
D_\nu=\sum_a\lfloor\nu(\{a\})\rfloor\delta_a,\qquad\nu^*=\nu-D_\nu.
\end{equation}
Here $\lfloor t\rfloor$ denotes the greatest integer not exceeding $t$.
The measure $D_\nu$ has locally finite support, since $D_\nu\leq\nu$ and each of its nonzero atoms has mass at least one. The remainder can have both a nonatomic part and atoms of mass strictly less than one.

\begin{lemma}\label{lem:residual}
Let $\sigma$ be a finite positive measure on $\C$. Define
\begin{equation}\label{eq:respotential}
v(z)=\int_\C\bigl(\log|z-\zeta|-\log\max(1,|\zeta|)\bigr)\dd\sigma(\zeta).
\end{equation}
Then $v$ is subharmonic and locally integrable on $\C$, its Riesz measure is $\sigma$, and
\begin{equation}\label{eq:Mv}
\inf_{z\in\C} M_1v(z)>-\infty.
\end{equation}
\end{lemma}
\begin{proof}
For $z$ in a fixed compact set and $|\zeta|$ sufficiently large, the kernel equals $\log|1-z/\zeta|$, which is uniformly $O(|z|/|\zeta|)$. On the other hand, for $\zeta$ in a bounded set, the kernel is locally integrable in $z$, and the integral of its absolute value over any compact set is bounded uniformly in $\zeta$. These estimates show that $v\in L^1_{\mathrm{loc}}(\C)$ and justify the use of Fubini's theorem below. To see upper semicontinuity, note that on $|z|\leq A$,
\[
\log|z-\zeta|-\log\max(1,|\zeta|)\leq\log(1+A).
\]
Applying Fatou's lemma to the nonnegative difference between $\log(1+A)$ and the logarithmic kernel along any sequence $z_j\to z$ shows that $v$ is upper semicontinuous. For any $\varphi\in C_c^\infty(\C)$, we have
\begin{align*}
\langle\Delta v,\varphi\rangle
&=\int_\C v(z)\Delta\varphi(z)\dd A(z)\\
&=\int_\C\int_\C
\bigl(\log|z-\zeta|-\log\max(1,|\zeta|)\bigr)\\
&\qquad\qquad\cdot\Delta\varphi(z)\dd A(z)\dd\sigma(\zeta)\\
&=2\pi\int_\C\varphi(\zeta)\dd\sigma(\zeta).
\end{align*}
The last equality uses $\int_\C\Delta\varphi\dd A=0$ and $\Delta_z\log|z-\zeta|=2\pi\delta_\zeta$. Hence $\Delta v=2\pi\sigma$. Integrating the submean inequality for $\log|z-\zeta|$ against $\sigma$, justified by the preceding local integrability bounds, shows that $v$ satisfies the submean inequality. Together with upper semicontinuity, this proves that $v$ is subharmonic. The distributional calculation gives $\nu_v=\sigma$. The function $v$ may take the value $-\infty$ at some points.

Let
\[
q(t)=\frac1\pi\int_{D(0,1)}\log|w-t|\dd A(w)
=\begin{cases}\log t,&t\geq1,\\(t^2-1)/2,&0\leq t\leq1.\end{cases}
\]
The preceding integrability estimates permit Fubini's theorem and give
\[
M_1v(z)=\int_\C K(z,\zeta)\dd\sigma(\zeta),
\]
where
\[
K(z,\zeta)=q(|z-\zeta|)-\log\max(1,|\zeta|).
\]
Since
\[
q(t)\geq\log\max(1,t)-\frac12,\qquad
\max(1,|\zeta|)\leq(1+|z|)\max(1,|z-\zeta|),
\]
we have
\begin{equation}\label{eq:Kbound}
K(z,\zeta)\geq-\log(1+|z|)-\frac12.
\end{equation}
If $\sigma(\C)=0$, the assertion is immediate. Otherwise choose $R\geq1$ such that $\sigma(D_R)>\sigma(\C)/2$. If $|z|\geq2R+2$ and $|\zeta|<R$, then $K(z,\zeta)\geq\log|z|-C_R$. On the complement use \eqref{eq:Kbound}. It follows that
\[
\begin{aligned}
M_1v(z)
&\geq\sigma(D_R)(\log|z|-C_R)\\
&\quad-\bigl(\sigma(\C)-\sigma(D_R)\bigr)\left(\log(1+|z|)+\frac12\right)\\
&\geq\bigl(2\sigma(D_R)-\sigma(\C)\bigr)\log|z|-C.
\end{aligned}
\]
On a bounded $z$-set, \eqref{eq:Kbound} already gives a uniform lower bound. This proves \eqref{eq:Mv}.
\end{proof}

\begin{proposition}[Holomorphic realization]\label{prop:realize}
Let $u\in\SH(\C)$ satisfy $\int_\C e^{2u}\dd\gamma<\infty$. Then there exists a nonzero $F\in\F$ such that $n_F\geq D_{\nu_u}$. In particular, $\nu_u\geq n_\Lambda$ implies $N_\Lambda\neq\{0\}$.
\end{proposition}
\begin{proof}
Suppose first that $\nu_u^*(\C)=\infty$. We use \cite[Theorem~2.6]{BorichevLeYoussfi} with the weight $\psi=2u$, which is subharmonic on $\C$ and hence locally integrable. Its distributional Laplacian $\mu_\psi=4\pi\nu_u$ is positive and locally finite, hence $\sigma$-finite. In the notation of that theorem,
\[
\begin{aligned}
(\mu_\psi)^d
&=4\pi\sum_a\left\lfloor\frac{\mu_\psi(\{a\})}{4\pi}\right\rfloor\delta_a
=4\pi D_{\nu_u},\\
(\mu_\psi)^c&=\mu_\psi-(\mu_\psi)^d=4\pi\nu_u^*.
\end{aligned}
\]
Since $(\mu_\psi)^c(\C)=4\pi\nu_u^*(\C)=\infty$, \cite[Theorem~2.6]{BorichevLeYoussfi} yields
\[
\dim\left\{g\in\Hol(\C):\int_\C|g|^2e^{-2u}\dd A<\infty\right\}=\infty.
\]
Choose a nonzero member $g$, so that
\[I:=\int_\C|g|^2e^{-2u}\dd A<\infty.\]
We first show that $n_g\geq D_{\nu_u}$. For any $a\in\C$, put $\alpha=\nu_u(\{a\})$. Separating the atom at $a$ in the local Riesz decomposition, we may write
\[u(z)=\alpha\log|z-a|+v_a(z)\]
near $a$, where $v_a$ is subharmonic and hence bounded above on a smaller disk. If $g$ had a zero of order $k<\lfloor\alpha\rfloor$ at $a$, then for some $\delta>0$,
\[
I\geq\delta\int_0^\varepsilon r^{2k-2\alpha+1}\dd r=\infty,
\]
since $k\leq\alpha-1$, contradicting $I<\infty$.

Both $\log|g|$ and $u$ belong to $L^1_{\mathrm{loc}}(\C)$, so
$\log(|g|^2e^{-2u})=2\log|g|-2u$ is locally integrable almost everywhere.
Values on the polar sets where an individual factor is infinite do not affect the following averages.
The submean inequality for $\log|g|$, followed by Jensen's inequality, gives
\begin{align*}
2\log|g(z)|&\leq M_1(2\log|g|)(z)\\
&=2M_1u(z)+M_1\log(|g|^2e^{-2u})(z)\\
&\leq2M_1u(z)+\log(I/\pi).
\end{align*}
Thus Lemma~\ref{lem:gaussavg} gives
\[
\norm{g}_{\F}^2\leq\frac I\pi\int e^{2M_1u}\dd\gamma
\leq\frac{eI}{\pi}\int e^{2u}\dd\gamma<\infty.
\]
We may take $F=g$.

Now suppose that $\nu_u^*(\C)<\infty$. Applying Lemma~\ref{lem:residual} to $\nu_u^*$, we obtain a subharmonic function $v$ with $\nu_v=\nu_u^*$. By the Weierstrass theorem there is a nonzero entire function $P$ such that $n_P=D_{\nu_u}$. The locally integrable difference $u-\log|P|-v$ has zero distributional Laplacian. By Weyl's lemma it agrees almost everywhere with a real harmonic function $H$ on $\C$. Since $\C$ is simply connected, $H=\Ree h$ for an entire function $h$. Set $F=Pe^h$. Since both $u$ and $\log|F|+v$ are subharmonic and agree almost everywhere, they agree everywhere:
\[u=\log|F|+v.\]
By \eqref{eq:Mv}, there is a constant $C_0$ such that $M_1v(z)\geq-C_0$ for all $z$. Consequently,
\[
\log|F(z)|\leq M_1\log|F|(z)=M_1u(z)-M_1v(z)\leq M_1u(z)+C_0.
\]
A second application of Lemma~\ref{lem:gaussavg} yields
\[
\norm{F}_{\F}^2\leq e^{2C_0+1}\int e^{2u}\dd\gamma<\infty.
\]
In this case $n_F=D_{\nu_u}$. Finally, $\nu_u(\{\lambda\})\geq m_\lambda$ implies $D_{\nu_u}(\{\lambda\})\geq m_\lambda$, proving the last assertion.
\end{proof}

\subsection{Proof of the general criterion}
\begin{proof}[Proof of Theorem~\ref{thm:unique}]
Suppose that $N_\Lambda\neq\{0\}$. Then there exists $f\in N_\Lambda$ with $f\not\equiv0$. If $f(0)=0$, let $k\geq1$ be the order of its zero at the origin and replace $f$ by $f(z)/z^k$. The quotient is entire after removal of the singularity at zero. It belongs to $\F$, since $z^{-k}$ is bounded for $|z|\geq1$ and the quotient is bounded on the closed unit disk. Since $0\notin\Lambda$, all prescribed vanishing conditions are preserved. Multiplying by a constant, we may therefore assume $f(0)=1$. The function $u=\log|f|$ belongs to $U_\Lambda$, and
\[
\int_\C e^{2u}\dd\gamma
=\int_\C|f|^2\dd\gamma
=\norm{f}_{\F}^2<\infty.
\]
Lemma~\ref{lem:entropy} gives $\E(\Lambda)<\infty$.

Conversely, if $\E(\Lambda)<\infty$, Proposition~\ref{prop:duality} supplies $u\in U_\Lambda$ with finite exponential integral. Proposition~\ref{prop:realize} gives a nonzero $F\in N_\Lambda$. This proves the equivalence.
\end{proof}

\begin{proposition}[Finite replacements]\label{rem:finite-stability}
Let $\Lambda$ be a locally finite uniqueness sequence with distinct terms, identified with its range. Suppose that $S\subset\Lambda$ and $T\subset\C\setminus(\Lambda\setminus S)$ are finite sets of distinct points such that $|T|\geq|S|$. Then
\[
\Lambda'=(\Lambda\setminus S)\cup T
\]
is again a uniqueness sequence. In particular, finite replacements with the same number of removed and inserted terms preserve uniqueness. If $\Lambda$ has critical asymptotic density, then $\Lambda'$ retains that density because finite modifications do not alter it.
\end{proposition}

\begin{proof}
Put $m=|S|$. If $0\ne f\in N_{(\Lambda\setminus S)\cup T}$, select distinct
$t_1,\ldots,t_m\in T$ and define
\[
g(z)=f(z)\frac{\prod_{s\in S}(z-s)}{\prod_{j=1}^{m}(z-t_j)}.
\]
The apparent poles are removable because $f(t_j)=0$, so $g$ is a nonzero entire function.
The quotient of the two monic degree-$m$ polynomials is bounded outside a disk; hence $g\in\F$.
At each point of the original sequence outside $S$, the denominator is nonzero and $f$ vanishes.
At each $s\in S$ not selected among the $t_j$, the numerator forces a zero.
If $s=t_j$, the two factors cancel and $f(s)=0$ because $s\in T$.
Thus $g$ vanishes on every term of the original sequence, contradicting uniqueness.
\end{proof}
The critical lattice example after Theorem~\ref{thm:angular} shows that the cardinality hypothesis cannot be dropped.

\section{Radial criteria and critical configurations}\label{sec:radial}

\subsection{The radial formula}
For $\Lambda$ as in Theorem~\ref{thm:unique}, let
\begin{equation}\label{eq:Phi}
\Phi_\Lambda(r)=\sum_{|\lambda|\leq r}m_\lambda\log\frac r{|\lambda|},
\qquad r>0.
\end{equation}
This function is continuous, finite on bounded intervals, and zero near the origin. Denote by $\E_{\mathrm{rad}}(\Lambda)$ the supremum defining $\E(\Lambda)$ when the test densities are restricted to radial densities.

\begin{proposition}[Radial formula]\label{prop:radial}
For a sequence $\Lambda$ as in Theorem~\ref{thm:unique},
\begin{equation}\label{eq:radial}
\E_{\mathrm{rad}}(\Lambda)=
\log\left(2\int_0^\infty re^{2\Phi_\Lambda(r)-r^2}\dd r\right).
\end{equation}
In particular, divergence of this integral implies that $\Lambda$ is a uniqueness sequence in $\F$.
\end{proposition}
\begin{proof}
If $\rho$ is radial, then $p(r)=2\pi r\rho(r)$ is a probability density on $(0,\infty)$. The circular mean identity gives
\[
\sum_{\lambda\in\operatorname{ran}\Lambda}m_\lambda V_\rho(\lambda)
=\int_0^\infty\Phi_\Lambda(r)p(r)\dd r,
\qquad
H_\gamma(\rho)=\int_0^\infty p(r)\log\frac{p(r)}{2re^{-r^2}}\dd r.
\]
Put $Z=2\int_0^\infty re^{2\Phi_\Lambda(r)-r^2}\dd r$. If $Z<\infty$, the density
\[
q(r)=Z^{-1}2re^{2\Phi_\Lambda(r)-r^2}
\]
satisfies
\[
2\int_0^\infty\Phi_\Lambda p\dd r
-\int_0^\infty p\log\frac p{2re^{-r^2}}\dd r
=\log Z-\int_0^\infty p\log\frac pq\dd r
\leq\log Z.
\]
This proves the upper bound.

For the reverse bound, fix $0<a<b<\infty$ and normalize
\[
p_{a,b}(r)=
\frac{2re^{2\Phi_\Lambda(r)-r^2}\ind_{[a,b]}(r)}
     {\displaystyle2\int_a^b te^{2\Phi_\Lambda(t)-t^2}\dd t}.
\]
Its variational value is $\log\left(2\int_a^b re^{2\Phi_\Lambda(r)-r^2}\dd r\right)$.
Let $\eta_\varepsilon$ be a smooth probability mollifier on $\R$, supported in $(-\varepsilon,\varepsilon)$, where $0<\varepsilon<a/2$, and set
\[
p_\varepsilon=p_{a,b}*\eta_\varepsilon,
\qquad
\rho_\varepsilon(z)=\frac{p_\varepsilon(|z|)}{2\pi|z|}.
\]
Here $\rho_\varepsilon(0)=0$; the density vanishes on a neighborhood of zero. Then $\rho_\varepsilon\in\J$. The functions $p_\varepsilon$ are uniformly bounded, converge almost everywhere to $p_{a,b}$, and have support in a fixed compact subinterval of $(0,\infty)$. Since $x\log x$ is continuous at zero with value zero, dominated convergence applies to both terms of the variational expression. Thus the values of $\rho_\varepsilon$ tend to the value of $p_{a,b}$. Letting $a\downarrow0$ and $b\uparrow\infty$ proves \eqref{eq:radial}, whether $Z$ is finite or infinite. Finally, $\E\geq\E_{\mathrm{rad}}$ and Theorem~\ref{thm:unique} gives the uniqueness assertion.
\end{proof}

Divergence implies uniqueness directly by Jensen's formula and Jensen's inequality (see \cite[proof of Theorem~2.1]{KononovaRandom} and \cite[Section~2.1, equation~(3)]{HouMa}). Formula~\eqref{eq:radial} gives the value of the radial supremum. The two configurations in \cite[Theorem~8.7]{HouMa}, with the common origin removed if present, have identical radial counts and hence the same radial supremum. This value is finite because one configuration is an exact zero sequence. The other is a uniqueness sequence, so its full supremum is infinite by Theorem~\ref{thm:unique}.

\subsection{Lattices below the critical area}
\begin{corollary}[Lattices of area less than $\pi$]\label{cor:subcritical}
Let $\Gamma=a+\omega_1\Z+\omega_2\Z$, where
$A=\im(\overline{\omega_1}\omega_2)\in(0,\pi)$. Regard $\Gamma$ as the sequence listing each lattice point once. Then $\Gamma$, and every finite modification of it, is a uniqueness sequence for $\F$.
\end{corollary}
\begin{proof}
Let $\Gamma'$ be a finite modification of $\Gamma$, and put $\Lambda=\Gamma'\setminus\{0\}$. Choose a bounded half-open fundamental parallelogram $Q$ of the period lattice and $D>0$ such that $Q\subset D(0,D)$. The cells $\gamma+Q$, $\gamma\in\Gamma$, partition the plane. For $r>D$, the union of the cells with $|\gamma|<r$ contains $D(0,r-D)$ and is contained in $D(0,r+D)$. Comparing areas gives
\[
\frac{\pi(r-D)^2}{A}
\leq\#\{\gamma\in\Gamma:|\gamma|<r\}
\leq\frac{\pi(r+D)^2}{A}.
\]
Finite modifications change the count by $O(1)$. Hence
\[
N_\Lambda(r):=\#\{\lambda\in\Lambda:|\lambda|<r\}
=\frac\pi A r^2+O(r),
\qquad
\Phi_\Lambda(r)=\int_0^r\frac{N_\Lambda(t)}t\dd t
=\frac\pi{2A}r^2+O(r).
\]
Since $A<\pi$, we have $2\Phi_\Lambda(r)-r^2\geq cr^2$ for all sufficiently large $r$ and some $c>0$. The radial integral diverges, so Proposition~\ref{prop:radial} gives uniqueness for $\Lambda$ and hence for $\Gamma'$.
\end{proof}

At $A=\pi$, this counting estimate leaves an $O(r)$ error. The critical case requires the sharper estimate proved below.

\subsection{Rotation orbits}\label{sec:orbits}
The radial condition is also necessary when the prescribed zeros have a generating function bounded by a constant multiple of $e^{\Phi_\Lambda(|z|)}$. A bounded overlap condition gives such a function explicitly. All vanishing conditions in this subsection are simple.

\begin{theorem}[Rotation orbits]\label{thm:orbits}
Let
\[
\mathcal O_j=\{r_je^{i(\alpha_j+2\pi k/N_j)}:0\leq k<N_j\},
\qquad r_j>0,\quad N_j\in\N,\quad\alpha_j\in\R,
\]
be pairwise disjoint rotation orbits. Assume that only finitely many orbits have radius at most any given $R<\infty$, and that
\begin{equation}\label{eq:overlap}
B:=\sup_{r>0}\sum_j\exp\bigl(-N_j|\log(r/r_j)|\bigr)<\infty.
\end{equation}
Let $\Lambda$ be the sequence obtained by listing each point of $\bigcup_j\mathcal O_j$ once. Its radial counting potential is
\[
\Phi_\Lambda(r)=\sum_{r_j\leq r}N_j\log(r/r_j).
\]
Then
\begin{equation}\label{eq:orbitcriterion}
N_\Lambda=\{0\}
\quad\Longleftrightarrow\quad
\int_1^\infty re^{2\Phi_\Lambda(r)-r^2}\dd r=\infty.
\end{equation}
If the integral is finite, the product in \eqref{eq:orbitproduct} belongs to $\F$ and has precisely the simple zeros of $\Lambda$. Thus the uniqueness property is independent of the angles $\alpha_j$, provided the orbits remain disjoint.
\end{theorem}
\begin{proof}
Set $a_j=r_je^{i\alpha_j}$ and consider
\begin{equation}\label{eq:orbitproduct}
P(z)=\prod_j\bigl(1-(z/a_j)^{N_j}\bigr).
\end{equation}
For each $R>0$, condition \eqref{eq:overlap} gives
$\sum_{r_j>R}(R/r_j)^{N_j}<\infty$. Hence the product converges locally uniformly, satisfies $P(0)=1$, and has exactly the stated simple zeros. For $|z|=r$,
\[
|1-(z/a_j)^{N_j}|
\leq\max\{1,(r/r_j)^{N_j}\}
\left(1+e^{-N_j|\log(r/r_j)|}\right).
\]
Taking products and using $\log(1+t)\leq t$ gives
\begin{equation}\label{eq:productbound}
|P(z)|\leq e^B e^{\Phi_\Lambda(|z|)}.
\end{equation}
Consequently, finiteness of the integral implies $P\in\F$ and nonuniqueness.

If the integral in \eqref{eq:orbitcriterion} diverges, Proposition~\ref{prop:radial} gives uniqueness.
\end{proof}

\begin{lemma}\label{lem:overlap}
Condition \eqref{eq:overlap} holds if $N_j\geq cr_j$ for some $c>0$ and all but finitely many $j$, and
\begin{equation}\label{eq:unitoverlap}
M:=\sup_{R\geq0}\#\{j:R\leq r_j<R+1\}<\infty.
\end{equation}
\end{lemma}
\begin{proof}
The exceptional indices contribute a uniformly bounded sum. For all other indices with $r_j<r/2$ or $r_j>2r$,
\[
e^{-N_j|\log(r/r_j)|}\leq e^{-cr_j\log2}.
\]
Their sum is at most $M\sum_{k\geq0}e^{-ck\log2}$, uniformly in $r$. For $r/2\leq r_j\leq2r$,
\[
N_j|\log(r/r_j)|
\geq cr_j\frac{|r-r_j|}{\max\{r,r_j\}}
\geq\frac c2|r-r_j|.
\]
Each shell $k\leq|r-r_j|<k+1$ is contained in two closed intervals of length at most one. Each such interval contains at most $2M$ radii, by splitting off its right endpoint in \eqref{eq:unitoverlap}. Thus the remaining sum is at most $4M\sum_{k\geq0}e^{-ck/2}$.
\end{proof}

For example, take two $n$-point orbits at radii $n$ and $n+n^{-3}$ for each $n\geq2$. The lemma applies, while aligned rotations give distances between corresponding points tending to zero. Uniform separation is therefore unnecessary. The theorem also allows disjoint orbits at the same radius, subject to \eqref{eq:overlap}.

\subsection{Critical lattice radii}\label{sec:angular}
The following estimate extends a square-lattice bound used in \cite[proof of Corollary~5.4]{HouMa}.

\begin{lemma}\label{lem:logdist}
Let $\Gamma$ be a full-rank lattice in $\C$, and write $d_\Gamma(z)=\dist(z,\Gamma)$. There is a constant $C_\Gamma<\infty$ such that
\begin{equation}\label{eq:logdist}
\left|\frac1{2\pi}\int_0^{2\pi}\log d_\Gamma(re^{it})\dd t\right|
\leq C_\Gamma,\qquad r\geq1.
\end{equation}
The logarithmic singularities where the circle meets the lattice are integrable.
\end{lemma}
\begin{proof}
Let $s>0$ be the minimum lattice separation, and choose
$0<\varepsilon_0<\min\{s/4,1/4\}$. There are $O(r)$ lattice points in
\begin{equation}\label{eq:annulus}
\{z:r-\varepsilon_0<|z|<r+\varepsilon_0\},\qquad r\geq1.
\end{equation}
Indeed, the disjoint disks $D(\gamma,s/3)$ centered at these points lie in the annulus with outer radius $r+\delta$ and inner radius $(r-\delta)_+$, where $\delta=\varepsilon_0+s/3$. Their total area is at most
\[
\pi\bigl((r+\delta)^2-(r-\delta)_+^2\bigr)\leq C_\Gamma r.
\]

Fix $0<\varepsilon\leq\varepsilon_0$. If $D(\gamma,\varepsilon)$ meets $|z|=r$, write $\gamma=be^{i\alpha}$ and choose the principal angular difference $|t-\alpha|\leq\pi$. On the intersecting arc,
\[
(r-b)^2+4rb\sin^2((t-\alpha)/2)<\varepsilon^2,
\qquad b>r-\varepsilon\geq3r/4.
\]
Since $|(t-\alpha)/2|\leq\pi/2$, the inequality
$|\sin x|\geq2|x|/\pi$ gives
\[
\left|\frac{t-\alpha}{2}\right|
\leq\frac{\pi\varepsilon}{4\sqrt{rb}}
\leq\frac{\pi\varepsilon}{2\sqrt3\,r}.
\]
The arc has angular length $O(\varepsilon/r)$. Summing over the $O(r)$ possible centers proves
\begin{equation}\label{eq:sublevel}
\frac1{2\pi}\bigl|\{t\in[0,2\pi):d_\Gamma(re^{it})<\varepsilon\}\bigr|
\leq C\varepsilon,
\qquad r\geq1,\quad0<\varepsilon\leq\varepsilon_0.
\end{equation}
For $Y_r=(\log(\varepsilon_0/d_\Gamma(re^{it})))_+$, the layer-cake formula and Tonelli's theorem imply
\[
\frac1{2\pi}\int_0^{2\pi}Y_r(t)\dd t
=\int_0^\infty\frac1{2\pi}
\bigl|\{t:d_\Gamma(re^{it})<\varepsilon_0e^{-v}\}\bigr|\dd v
\leq C\varepsilon_0.
\]
Since $\log d_\Gamma\geq\log\varepsilon_0-Y_r$, this is the required lower bound. The finite covering radius of $\Gamma$ gives the upper bound. These estimates also prove integrability at the intersections.
\end{proof}

\begin{proposition}\label{prop:criticrad}
Let $\Gamma=\omega_1\Z+\omega_2\Z$ have fundamental area $\pi$. Then
\begin{equation}\label{eq:criticPhi}
\Phi_{\Gamma\setminus\{0\}}(r)
:=\sum_{\substack{\gamma\in\Gamma\\0<|\gamma|\leq r}}
\log\frac r{|\gamma|}
=\frac{r^2}{2}-\log r+O(1),\qquad r\geq1.
\end{equation}
The error is bounded uniformly in $r$ for the fixed lattice.
\end{proposition}
\begin{proof}
Choose positively oriented generators, so that $\im(\overline{\omega_1}\omega_2)=\pi$. Let $\sigma_\Gamma$ be the Weierstrass sigma function of $\Gamma$. In full-period notation, quasiperiodicity and the Legendre relation give
\[
\sigma_\Gamma(z+\omega_j)
=-e^{\eta_j(z+\omega_j/2)}\sigma_\Gamma(z),
\qquad
\eta_1\omega_2-\eta_2\omega_1=2\pi i;
\]
see \cite[Proposition~2.5(ii)--(iii)]{GrochenigLyubarskii}. Since
$\overline{\omega_1}\omega_2-\overline{\omega_2}\omega_1=2\pi i$,
there is $c\in\C$ such that
$\eta_j-2c\omega_j=\overline{\omega_j}$ for $j=1,2$.
Thus $S_\Gamma(z)=e^{-cz^2}\sigma_\Gamma(z)$ satisfies
\[
S_\Gamma(z+\omega_j)
=-\exp\!\left(\overline{\omega_j}(z+\omega_j/2)\right)S_\Gamma(z).
\]
It follows that $|S_\Gamma(z)|e^{-|z|^2/2}$ is $\Gamma$-periodic. Its zeros are simple and are precisely the points of $\Gamma$. Dividing by $d_\Gamma$ gives a positive continuous periodic function, after extension at the lattice points. Hence
\[
|S_\Gamma(z)|e^{-|z|^2/2}\asymp d_\Gamma(z).
\]
Normalize $S_\Gamma'(0)=1$. Taking logarithms and circular averages, justified by Lemma~\ref{lem:logdist}, yields
\[
\frac1{2\pi}\int_0^{2\pi}\log|S_\Gamma(re^{it})|\dd t
=\frac{r^2}{2}+O(1).
\]
The function $S_\Gamma(z)/z$ equals one at zero and has precisely the zeros $\Gamma\setminus\{0\}$. Jensen's formula therefore gives
\[
\Phi_{\Gamma\setminus\{0\}}(r)
=\frac1{2\pi}\int_0^{2\pi}\log|S_\Gamma(re^{it})|\dd t-\log r,
\]
which is \eqref{eq:criticPhi}. The formula remains valid on circles containing zeros: the logarithmic singularities are integrable, and boundary zeros contribute zero to the sum.
\end{proof}

\begin{theorem}[Angular rearrangements]\label{thm:angular}
Let $\Gamma$ be a full-rank lattice containing the origin, with fundamental area $\pi$. For each $\gamma\in\Gamma\setminus\{0\}$ choose $\lambda_\gamma$ with $|\lambda_\gamma|=|\gamma|$, and assume that the chosen points are pairwise distinct. Then
\[
\Lambda^\circ=(\lambda_\gamma)_{\gamma\in\Gamma\setminus\{0\}}
\]
is a uniqueness sequence. The sequence obtained by adjoining the origin to $\Lambda^\circ$ is also a uniqueness sequence.
\end{theorem}
\begin{proof}
The sequence $\Lambda^\circ$ is locally finite and avoids zero. Since $|\lambda_\gamma|=|\gamma|$ for each label, the two radial counting potentials agree. Proposition~\ref{prop:criticrad} gives
\begin{equation}\label{eq:oneoverr}
\frac c r\leq re^{2\Phi_{\Lambda^\circ}(r)-r^2}\leq\frac C r,
\qquad r\geq1,
\end{equation}
for constants $0<c\leq C<\infty$. The radial integral diverges, and Proposition~\ref{prop:radial} proves that $\Lambda^\circ$ is a uniqueness sequence. Adding the origin preserves this property.
\end{proof}

The endpoint is sharp for the unrearranged lattice. Let $S_\Gamma$ be the modified sigma function constructed in the proof of Proposition~\ref{prop:criticrad}. If $a\in\Gamma\setminus\{0\}$, then
\[
F_a(z)=\frac{S_\Gamma(z)}{z(z-a)}
\]
is entire, with simple zeros precisely at $\Gamma\setminus\{0,a\}$. The sigma estimate gives
$|F_a(z)|e^{-|z|^2/2}=O(|z|^{-2})$ at infinity, so $F_a\in\F$.
Thus $\Gamma\setminus\{0\}$ is a uniqueness sequence, whereas deleting one further point leaves a sequence that is not a uniqueness sequence. This also illustrates the role of the cardinality condition in Proposition~\ref{rem:finite-stability}.

\subsection{Independent rotations of square-lattice radii}
Chistyakov, Lyubarskii, and Pastur \cite{ChistyakovLyubarskiiPastur} studied random completeness in the Bargmann--Fock space in a different random model. Kononova \cite[Introduction, item~2(c)]{KononovaRandom} considered the following model. Let $\theta_{m,n}$, $(m,n)\in\Z^2\setminus\{(0,0)\}$, be independent and uniform on $[0,2\pi)$, and put
\begin{equation}\label{eq:randommodel}
\begin{aligned}
\lambda_{m,n}^{(a)}&=\sqrt{a\pi(m^2+n^2)}\,e^{i\theta_{m,n}},\\
\Lambda_a&=\bigl(0,\ (\lambda_{m,n}^{(a)})_{(m,n)\in\Z^2\setminus\{(0,0)\}}\bigr),\qquad a>0.
\end{aligned}
\end{equation}
The model gives exact zero sequences almost surely for $a>1$ and uniqueness sequences for $a<1$. Kononova asked what happens at $a=1$ \cite[Introduction, item~2(c)]{KononovaRandom}; see also \cite[Section~9]{KononovaI}. Theorem~\ref{thm:angular} settles the critical case.

\begin{corollary}\label{cor:random}
For the model \eqref{eq:randommodel} at $a=1$,
\begin{equation}\label{eq:almostsure}
\mathbb P\bigl(\{\omega:\ \forall F\in\F,\ 
F|_{\Lambda_1(\omega)}=0\Longrightarrow F\equiv0\}\bigr)=1.
\end{equation}
The same holds with the origin omitted. Equivalently, the associated coherent states are complete almost surely.
\end{corollary}
\begin{proof}
Different radii cannot collide. For two distinct labels of the same radius, independence and uniformity give collision probability zero. The labels have only countably many pairs, so on one event of probability one all nonzero points are pairwise distinct. Apply Theorem~\ref{thm:angular} on this event to $\Gamma=\sqrt\pi(\Z+i\Z)$.
\end{proof}

The coherent state at $\lambda$ is the normalized reproducing kernel
$k_\lambda(z)=e^{z\overline\lambda-|\lambda|^2/2}$. Orthogonality to $k_\lambda$ is equivalent to vanishing at $\lambda$. The completeness event is measurable. It can be written in terms of approximation of a countable dense set in $\F$ by finite linear combinations of the kernels with coefficients in $\Q+i\Q$, and $\lambda\mapsto k_\lambda$ is continuous in the norm of $\F$.

\section{Countable tests and finite geometric approximation}\label{sec:polynomials}

\subsection{Polynomial Jensen tests}
Let $m$ denote normalized arclength measure on the unit circle $\T$. Define the countable polynomial class
\[
\mathcal P_{\Q}=\left\{P(\zeta)=\sum_{k=1}^d a_k\zeta^k:
\ d\in\N,\ a_k\in\Q+i\Q\right\}.
\]
The zero polynomial is allowed. For a polynomial $P$, let $P_*m$ denote the pushforward of $m$ under $P$, so
\[
(P_*m)(E)=m(P^{-1}(E))=m\{\zeta\in\T:P(\zeta)\in E\}.
\]
For $P\in\mathcal P_{\Q}$ and $\delta\in\Q_{>0}$, put
\begin{equation}\label{eq:diskdensity}
\begin{gathered}
\mu_P=P_*m,\qquad \kappa_\delta(z)=\frac{\ind_{\{|z|<\delta\}}}{\pi\delta^2},\\
\rho_{P,\delta}(z)=(\mu_P*\kappa_\delta)(z)
=\frac{m\{\zeta\in\T:|z-P(\zeta)|<\delta\}}{\pi\delta^2}.
\end{gathered}
\end{equation}
The disk average of the logarithmic kernel is given by
\begin{equation}\label{eq:ell}
\ell_\delta(s)=\begin{cases}
\log\delta-\frac12+\dfrac{s^2}{2\delta^2},&0\leq s\leq\delta,\\
\log s,&s\geq\delta.
\end{cases}
\end{equation}
\begin{lemma}\label{lem:polyJ}
The density $\rho_{P,\delta}$ is a compactly supported probability density, bounded by $M_\delta=(\pi\delta^2)^{-1}$, and is a Jensen density at zero. It has finite Gaussian relative entropy. With
\[B=\delta+\sum_{k=1}^d|a_k|,\]
its potential is zero for $|w|>B$, and for $w\neq0$ it is given by
\begin{equation}\label{eq:polyV}
V_{P,\delta}(w)=\int_\T\ell_\delta(|w-P(\zeta)|)\dd m(\zeta)-\log|w|\geq0.
\end{equation}
Its entropy satisfies
\begin{equation}\label{eq:polyH}
H_\gamma(\rho_{P,\delta})
=\log\pi+\sum_{k=1}^d|a_k|^2+\frac{\delta^2}{2}
+\int_\C\rho_{P,\delta}\log\rho_{P,\delta}\dd A.
\end{equation}
\end{lemma}
\begin{proof}
The mass, bound, and support follow directly from \eqref{eq:diskdensity}. If $u$ is subharmonic, then $u\circ P$ is subharmonic on the unit disk, and
\[u(0)=u(P(0))\leq\int_\T u(P(\zeta))\dd m(\zeta)
\leq\int_\C u(z)\rho_{P,\delta}(z)\dd A(z).
\]
The second inequality follows from the disk submean inequality at each $P(\zeta)$. In particular, applying these inequalities to $u(z)=\log|w-z|$ proves that $V_{P,\delta}(w)\geq0$ for $w\neq0$.

The average of $\log|w-z|$ over $|z|<\delta$, divided by $\pi\delta^2$, is $\ell_\delta(|w|)$; this follows by integrating the circular mean $\log\max\{|w|,r\}$ against $2r\dd r/\delta^2$. Fubini's theorem gives \eqref{eq:polyV}. For $|w|>B$, disk averaging leaves the logarithm unchanged and $w-P$ has no zero in the closed unit disk. Its logarithmic mean is $\log|w|$, so $V_{P,\delta}(w)=0$.

Finally, the disk has mean zero and second moment $\delta^2/2$. Parseval's identity yields
\[
\int_\C|z|^2\rho_{P,\delta}(z)\dd A(z)
=\sum_{k=1}^d|a_k|^2+\frac{\delta^2}{2}.
\]
This proves \eqref{eq:polyH}. Boundedness and compact support, together with continuity of $x\log x$ at zero, prove finiteness.
\end{proof}

We justify the use of these densities in the supremum by smoothing. They provide nonradial tests. For $P(\zeta)=\zeta+\zeta^2/2$, direct integration gives
\[
\int_\C z|z|^2\rho_{P,\delta}(z)\dd A(z)=\frac12,
\]
whereas this mixed moment vanishes for every radial density. The disk convolution does not change the moment because $\int P\dd m=0$.

There is a second geometric form of the potential. For $P\not\equiv0$, define, counting roots with multiplicity,
\[
Z_P(w)=\sum_{\substack{P(\zeta)=w\\|\zeta|<1}}\log\frac1{|\zeta|},\qquad w\neq0,
\]
and put $Z_0=0$. Jensen's formula for $P-w$ gives
\begin{equation}\label{eq:polyalternate}
V_{P,\delta}(w)=\ell_\delta(|w|)-\log|w|
+\frac1{\pi\delta^2}\int_{|t|<\delta}Z_P(w-t)\dd A(t).
\end{equation}
The singular value at $w-t=0$ is interpreted through the locally integrable logarithmic potential. Formula \eqref{eq:polyV} is often simpler for evaluation because its kernel is bounded near zero.

For a concrete directional test, take $P(\zeta)=\zeta+\zeta^2/2$ and $\delta=1/4$. Then
\[
V_{P,1/4}(1)=\log\frac{1+\sqrt3}{2}>0,
\qquad V_{P,1/4}(-1)=0.
\]
Indeed, the distances from $1$ and $-1$ to $P(\T)$ are $1/2$ and $1/(2\sqrt2)$, respectively. The logarithmic kernel is harmonic in the averaging disks, since they avoid $P(\T)$ and the origin. Thus the smoothed and unsmoothed potentials agree at $1$ and $-1$. The roots of $P(\zeta)=1$ are $-1\pm\sqrt3$, whereas those of $P(\zeta)=-1$ are $-1\pm i$. Thus the test gives different weights to two points of the same modulus.

\subsection{An exact countable criterion}
For the rest of this section, $\Lambda$ is a locally finite sequence avoiding the origin, with repetitions allowed as in Theorem~\ref{thm:unique}. Set
\begin{equation}\label{eq:polyscore}
S_\Lambda(P,\delta)=2\sum_{|\lambda|\leq B}m_\lambda V_{P,\delta}(\lambda)
-H_\gamma(\rho_{P,\delta}),\qquad B=\delta+\sum_{k=1}^d|a_k|.
\end{equation}
\begin{theorem}[Countable criterion]\label{thm:polycriterion}
For a sequence $\Lambda$ as in Theorem~\ref{thm:unique},
\begin{equation}\label{eq:polysup}
\E(\Lambda)=\sup_{\substack{P\in\mathcal P_\Q\\\delta\in\Q_{>0}}}S_\Lambda(P,\delta).
\end{equation}
Consequently, $N_\Lambda=\{0\}$ if and only if these explicit scores are unbounded above.
\end{theorem}

We use the following consequence of Poletsky's representation theorem.
\begin{lemma}\label{lem:polyapprox}
For every $\rho\in\J$, there are $P_j\in\mathcal P_{\Q}$ such that $\mu_{P_j}$ converges weakly to $\rho\dd A$ and
\[\sup_j\max_{|\zeta|\leq1}|P_j(\zeta)|<\infty.\]
\end{lemma}
\begin{proof}
Choose $0<R<S$ so that $\supp\rho\subset D(0,R)$, and put $D=D(0,S)$.
We first check the Jensen condition on this bounded domain.
If $u\in\SH(D)$ and $u(0)>-\infty$, choose $R<S'<S$ and apply
Lemma~\ref{lem:jensen} on $D(0,S')$ to obtain
\[
u(0)\leq\int_\C u\rho\dd A.
\]
The inequality is immediate when $u(0)=-\infty$.
Thus $\rho\dd A$ is a Jensen measure on $D$, with compact support in $D$.
Applying the inequality to both signs of the harmonic coordinate functions
also gives $\int_\C z\rho(z)\dd A(z)=0$.

Poletsky's theorem \cite[Theorem~2.1]{Poletsky} supplies holomorphic maps
$f_j$ from the unit disk into $D$ whose boundary measures converge weakly
to $\rho\dd A$. Let $f_j^*$ denote their radial boundary values.
Their centers satisfy
\[
f_j(0)=\int_\T f_j^*\dd m\longrightarrow0.
\]
Since the maps are bounded by $S$, radial convergence and dominated
convergence permit a choice of $0<r_j<1$ with
\[
\|f_j(r_j\cdot)-f_j^*\|_{L^1(m)}<j^{-1}.
\]
Choose a Taylor polynomial $Q_j$ of $f_j(r_j\cdot)$ such that
\[
\|Q_j-f_j(r_j\cdot)\|_{L^\infty(\{|\zeta|\leq1\})}<j^{-1},
\qquad Q_j(0)=f_j(0).
\]
Approximating its nonconstant coefficients by rational complex numbers
gives $P_j\in\mathcal P_\Q$ with
\[
\|P_j-(Q_j-Q_j(0))\|_{L^\infty(\{|\zeta|\leq1\})}<j^{-1}.
\]
Consequently,
\[
\|P_j-f_j^*\|_{L^1(m)}\leq |f_j(0)|+3j^{-1}\longrightarrow0,
\qquad
\|P_j\|_{L^\infty(\{|\zeta|\leq1\})}\leq 2S+2.
\]
Uniform continuity of continuous test functions on this fixed compact disk
now shows that $\mu_{P_j}$ and the boundary measures of $f_j$ have the same
weak limit. This proves the lemma.
\end{proof}

\begin{proof}[Proof of Theorem~\ref{thm:polycriterion}]
First let $\sigma$ be any bounded, compactly supported Jensen probability density. Convolve it with a smooth radial probability density $\eta_\varepsilon$, supported in $|z|<\varepsilon$. Lemma~\ref{lem:smoothing} gives
\[
\sigma*\eta_\varepsilon\in\J,\quad V_{\sigma*\eta_\varepsilon}\geq V_\sigma,\quad
H_\gamma(\sigma*\eta_\varepsilon)\leq H_\gamma(\sigma)+\varepsilon^2.
\]
Thus its unsmoothed score is at most $\E(\Lambda)+\varepsilon^2$. Letting $\varepsilon\downarrow0$, and then taking $\sigma=\rho_{P,\delta}$, proves that the right side of \eqref{eq:polysup} is at most the left side.

Conversely, fix $\rho\in\J$ and the sequence from Lemma~\ref{lem:polyapprox}. Fix a positive rational $\delta$. Every circle has zero $\rho\dd A$-measure, so weak convergence and \eqref{eq:diskdensity} imply, for every $z\in\C$,
\[
\rho_{P_j,\delta}(z)\longrightarrow\sigma_\delta(z):=(\rho*\kappa_\delta)(z).
\]
These densities have a common compact support and a common bound $M_\delta$. Dominated convergence, applied to $x\log x$ and the Gaussian moment, gives
\[H_\gamma(\rho_{P_j,\delta})\longrightarrow H_\gamma(\sigma_\delta).\]
The kernel $z\mapsto\ell_\delta(|w-z|)$ is continuous and bounded on the common support. Hence the potentials also converge at each $w\neq0$. They all vanish outside one fixed disk, since the polynomial images are uniformly bounded. Only finitely many points of $\Lambda$ therefore occur in the scores, and those scores converge to the score of $\sigma_\delta$.

Finally let rational $\delta\downarrow0$. The smooth densities $\rho*\kappa_\delta$ converge uniformly to $\rho$, with a common compact support. Their entropies converge, and their potentials converge by local integrability of the logarithmic kernel. The supremum on the right of \eqref{eq:polysup} therefore dominates the score of every $\rho\in\J$. This proves equality. The last assertion follows from Theorem~\ref{thm:unique}.
\end{proof}

The order of approximation in this proof is essential: first fix the smoothing radius and pass to the disk approximation, then let the radius decrease. Weak convergence alone would not justify entropy convergence for the unsmoothed boundary measures.
\begin{corollary}\label{cor:smoothtests}
Fix a smooth radial probability density $\eta$ supported in the unit disk, and put $\eta_n(z)=n^2\eta(nz)$. The countable family
\[
\J^*=\{\rho_{P,\delta}*\eta_n:P\in\mathcal P_\Q,\ \delta\in\Q_{>0},\ n\in\N\}\subset\J
\]
is independent of $\Lambda$, and the supremum in \eqref{eq:E} is unchanged when $\J$ is replaced by $\J^*$.
\end{corollary}
\begin{proof}
Lemma~\ref{lem:smoothing} puts every member in $\J$ and makes its score at least $S_\Lambda(P,\delta)-n^{-2}$. Now use Theorem~\ref{thm:polycriterion} and let $n\to\infty$.
\end{proof}

\subsection{Finite approximations and error bounds}
We begin with an entropy estimate.
\begin{lemma}\label{lem:enterror}
Let $f,g$ be probability densities, both bounded by $M>0$ and supported in a set of area $A<\infty$. Put $t=\|f-g\|_{L^1}$. For $t>0$,
\begin{equation}\label{eq:enterror}
\left|\int_\C f\log f\dd A-\int_\C g\log g\dd A\right|
\leq t\left(1+|\log M|+\log\frac{MA}{t}\right).
\end{equation}
The right-hand side is defined to be zero when $t=0$.
\end{lemma}
\begin{proof}
For $0\leq x,y\leq M$, put $h=|x-y|$. Scaling to $[0,1]$ gives
\[
|x\log x-y\log y|\leq h\left(1+|\log M|+\log\frac{M}{h}\right),
\]
with the usual zero convention. Indeed, for $0\leq v\leq v+a\leq1$, the difference of $-s\log s$ at $v+a$ and $v$ decreases with $v$; its positive endpoint is $-a\log a$, and the magnitude of its negative endpoint is at most $a$. Integrate the displayed bound and use concavity of $h\mapsto h\log(M/h)$ on a set of area $A$. Here $t\leq MA$, so $\log(MA/t)\geq0$.
\end{proof}

Fix $P(\zeta)=\sum_{k=1}^d a_k\zeta^k$, $\delta>0$, and an integer $K>d$. Put
\begin{equation}\label{eq:finitecenters}
\begin{gathered}
\theta_\ell=\frac{2\pi(\ell+1/2)}K,\qquad
c_\ell=P(e^{i\theta_\ell}),\\
I_\ell=\{e^{i\theta}:2\ell\pi/K\leq\theta\leq2(\ell+1)\pi/K\},\quad0\leq\ell<K,\\
D_P=\sum_{k=1}^d k|a_k|,\qquad
B=\delta+\sum_{k=1}^d|a_k|,\qquad
\Delta_K=\frac{\pi D_P}K.
\end{gathered}
\end{equation}
For $1\leq j\leq K$, let $A_j$ be the planar area covered by exactly $j$ of the $K$ disks $D(c_\ell,\delta)$, counting coincident disks with multiplicities. Define
\begin{equation}\label{eq:finitescorecomponents}
\begin{gathered}
\widehat\rho_K(z)=\frac1K\sum_{\ell=0}^{K-1}\kappa_\delta(z-c_\ell),\\
\widehat h_K=\sum_{j=1}^K A_j\frac{j}{K\pi\delta^2}\log\frac{j}{K\pi\delta^2}
=\int_\C\widehat\rho_K\log\widehat\rho_K\dd A,\\
\widehat V_K(w)=\frac1K\sum_{\ell=0}^{K-1}\ell_\delta(|w-c_\ell|)-\log|w|.
\end{gathered}
\end{equation}
The corresponding finite score is
\begin{equation}\label{eq:finitescore}
\widehat S_\Lambda(P,\delta,K)
=2\sum_{|\lambda|\leq B}m_\lambda\widehat V_K(\lambda)
-\log\pi-\sum_{k=1}^d|a_k|^2-\frac{\delta^2}2-\widehat h_K.
\end{equation}
The restriction $|\lambda|\leq B$ is part of the definition. No values of $\widehat V_K$ outside that disk are used.

Write $M=(\pi\delta^2)^{-1}$, $A=\pi B^2$, and set
\begin{equation}\label{eq:finiteerror}
\begin{gathered}
t_K=\min\{2,2\Delta_K/\delta,MA\},\\
\Omega_{M,A}(t)=t\left(1+|\log M|+\log\frac{MA}{t}\right),\quad0<t\leq MA,
\qquad \Omega_{M,A}(0)=0,\\
n_\Lambda(\overline D_B)=\sum_{|\lambda|\leq B}m_\lambda,\qquad
\varepsilon_\Lambda(P,\delta,K)
=\frac{2n_\Lambda(\overline D_B)\Delta_K}{\delta}+\Omega_{M,A}(t_K).
\end{gathered}
\end{equation}
\begin{theorem}[Certified finite tests]\label{thm:finitegeometry}
For $P\in\mathcal P_{\Q}$, $\delta\in\Q_{>0}$, and $K>d$,
\begin{equation}\label{eq:finiteerrorbound}
\left|S_\Lambda(P,\delta)-\widehat S_\Lambda(P,\delta,K)\right|
\leq\varepsilon_\Lambda(P,\delta,K).
\end{equation}
In particular, for $\Lambda$ as in Theorem~\ref{thm:unique},
\begin{equation}\label{eq:finitemax}
\E(\Lambda)=\sup_{\substack{P\in\mathcal P_\Q,\ \delta\in\Q_{>0}\\K\in\N,\ K>d}}
\left\{\widehat S_\Lambda(P,\delta,K)-\varepsilon_\Lambda(P,\delta,K)\right\}.
\end{equation}
We have $N_\Lambda=\{0\}$ if and only if these certified finite geometric scores are unbounded above.
\end{theorem}
\begin{proof}
The curve satisfies
\[
|P(e^{i\theta})-P(e^{i\phi})|\leq D_P|\theta-\phi|
\]
on each parametrizing interval. Couple each point of the unit circle to the midpoint of its interval of length $2\pi/K$. The displacement of its image is at most $\Delta_K$. The radial function $z\mapsto\ell_\delta(|z|)$ is $\delta^{-1}$-Lipschitz: its radial derivative is $s/\delta^2$ for $s\leq\delta$, and $1/s$ for $s\geq\delta$. Consequently,
\begin{equation}\label{eq:finiteVerror}
|V_{P,\delta}(w)-\widehat V_K(w)|\leq\Delta_K/\delta,\qquad w\neq0.
\end{equation}
For a translate by $h$, the normalized disk density satisfies
\[
\|\kappa_\delta(\cdot-h)-\kappa_\delta\|_{L^1}
\leq\min\{2,2|h|/\delta\}.
\]
Indeed, slicing two horizontally translated disks into intervals bounds the area of their symmetric difference by $4\delta|h|$. Dividing by $\pi\delta^2$ gives the translation bound, and the triangle inequality gives the bound by $2$. Since $|P(\zeta)-c_\ell|\leq\Delta_K$ for $\zeta\in I_\ell$, the coupling gives
\[
\begin{aligned}
\|\rho_{P,\delta}-\widehat\rho_K\|_{L^1}
&\leq\sum_{\ell=0}^{K-1}\int_{I_\ell}\!\int_\C
|\kappa_\delta(z-P(\zeta))-\kappa_\delta(z-c_\ell)|\dd A(z)\dd m(\zeta)\\
&\leq\min\{2,2\Delta_K/\delta\}.
\end{aligned}
\]
Both densities are bounded by $M$ and supported in $|z|\leq B$; their $L^1$ difference is also at most $MA$. The function $\Omega_{M,A}$ is nondecreasing on $[0,MA]$. Lemma~\ref{lem:enterror} therefore bounds the absolute difference of their integrals of $\rho\log\rho$ by $\Omega_{M,A}(t_K)$. Moreover, $K>d$ makes the midpoint quadrature exact for $|P(e^{i\theta})|^2$, so their Gaussian second moments agree with the value in \eqref{eq:polyH}. Together with \eqref{eq:finiteVerror}, these facts prove \eqref{eq:finiteerrorbound}.

By the absolute error estimate, each corrected score on the right of \eqref{eq:finitemax} is at most $S_\Lambda(P,\delta)$, hence at most $\E(\Lambda)$. For fixed $P,\delta$, local finiteness gives $n_\Lambda(\overline D_B)<\infty$, and $\varepsilon_\Lambda(P,\delta,K)\to0$ as $K\to\infty$. The same error estimate implies $\widehat S_\Lambda(P,\delta,K)\to S_\Lambda(P,\delta)$, so the supremum over $K$ of the corrected scores equals $S_\Lambda(P,\delta)$. Now apply Theorem~\ref{thm:polycriterion}.
\end{proof}

\subsection{Finite evaluation and increasing lower bounds}
Every quantity in \eqref{eq:finitemax} uses finitely many points of $\Lambda$ with their repetitions, distances to the explicitly specified centers $c_\ell$, and areas of a finite circle arrangement. These areas are obtained by intersecting the circles and summing oriented boundary arcs. For an arc of the circle centered at $x+iy$, with radius $\delta$ and angles $\alpha\leq\theta\leq\beta$, its contribution to Green's area formula is
\[
\frac12\left\{\delta^2(\beta-\alpha)+\delta x(\sin\beta-\sin\alpha)+\delta y(\cos\alpha-\cos\beta)\right\}.
\]
Compute first the superlevel areas
\[
C_j=A\bigl(\{z:\textstyle\sum_{\ell=0}^{K-1}
\ind_{D(c_\ell,\delta)}(z)\geq j\}\bigr),\qquad 1\leq j\leq K,
\]
and set $C_{K+1}=0$. Then $A_j=C_j-C_{j+1}$.
Group coincident circles, keeping their multiplicities, and subdivide
each remaining circle at its finitely many intersection points.
The coverage numbers on the two sides of each open arc are constant.
Give the arc the orientation for which the inside of its disk lies to the
left. If its two coverage numbers are $q_+$ on the left and $q_-$ on the
right, its Green integral contributes to $C_j$ exactly when
$q_-<j\leq q_+$. Summing these contributions computes all $C_j$ and hence
all $A_j$. Full circles with no intersection points are handled in the
same way; isolated tangencies have zero area. This procedure automatically
assigns the correct orientation along holes.

\begin{remark}
The discrete density $\widehat\rho_K$ is an evaluation device and need not be a Jensen density. In particular, $\widehat V_K$ need not be nonnegative or vanish outside $|w|\leq B$. The truncation in \eqref{eq:finitescore} and the error subtraction in \eqref{eq:finitemax} are essential. Replacing the exact curve measure by finitely many atoms without this control would not give a valid uniqueness criterion.
\end{remark}

For a rational test, the centers and the radius are algebraic. Thus coincidences, tangencies, and the cyclic order of circle intersections can be determined by exact algebraic comparisons. The arc integrals can then be enclosed to any prescribed accuracy.
To apply the test to an infinite sequence, however, one also needs an
effective description of that sequence. We assume a bounded-window
procedure which, for each test radius $B$, returns the complete finite list
in $|z|\leq B$, with exact multiplicities and certified rational enclosures
of every coordinate to any requested precision. Boundary membership is part
of this procedure. Local finiteness by itself does not provide such an input.

\begin{proposition}[Effective evaluation]\label{prop:effective}
Assume that the prescribed sequence is supplied by the bounded-window procedure above. For each $P\in\mathcal P_\Q$ and $\delta\in\Q_{>0}$, the score $S_\Lambda(P,\delta)$ is computable to any prescribed accuracy. Moreover, one can construct a nondecreasing sequence of nonnegative rational numbers $b_n$ such that
\[
b_n\uparrow\E(\Lambda),
\]
with $+\infty$ allowed as the limit.
\end{proposition}
\begin{proof}
For a single score, use the following finite procedure.
\begin{enumerate}
\item Choose $P\in\mathcal P_{\Q}$, $\delta\in\Q_{>0}$, and a tolerance
$\tau\in\Q_{>0}$. Obtain the complete list in $|z|\leq B$, and compute
$n_\Lambda(\overline D_B)$ and the constants in \eqref{eq:finiteerror}.
\item Increase $K>d$ until a certified upper bound for
$\varepsilon_\Lambda(P,\delta,K)$ is less than $\tau/4$.
This search terminates for each fixed pair $(P,\delta)$, since the error
tends to zero as $K\to\infty$.
\item Form the $K$ disks and evaluate their coverage areas by the boundary-arc
formula above. Coincident disks are retained with their multiplicities;
tangencies and empty coverage levels are handled in the same arrangement.
Use interval bounds for the areas and for every logarithm and distance in
\eqref{eq:finitescore}. If $|A_j-\widetilde A_j|\leq\eta_j$, the contribution
of these area errors to the entropy error is at most
\[
\sum_{j=1}^{K}\left|
\frac{j}{K\pi\delta^2}\log\frac{j}{K\pi\delta^2}
\right|\eta_j.
\]
Refine these bounds until the resulting interval for
$\widehat S_\Lambda(P,\delta,K)$ has width less than $\tau/4$.
\item Write $[s^-,s^+]$ for that interval, and let $e^+$ be the
certified upper bound for the quadrature error. Output
\[
L=s^- -e^+.
\]
Then
\[
S_\Lambda(P,\delta)-\tau<L\leq S_\Lambda(P,\delta)\leq\E(\Lambda),
\]
because $S_\Lambda(P,\delta)-L\leq(s^+-s^-)+2e^+<3\tau/4$.
All coordinate, area, and elementary-function errors are included in this bound.
Thus $[L,L+\tau]$ is a certified interval for the exact score.
\end{enumerate}

Enumerate the scores $S_\Lambda(P,\delta)$ in \eqref{eq:polysup} as $T_1,T_2,\ldots$, and let $\ell_{j,n}$ be the certified rational lower approximations from the preceding procedure, chosen so that
\[
T_j-2^{-n}\leq\ell_{j,n}\leq T_j.
\]
Set $b_0=0$ and
\[
b_n=\max\{b_{n-1},0,\ell_{1,n},\ldots,\ell_{n,n}\}.
\]
Since $\E(\Lambda)\geq0$, Theorem~\ref{thm:polycriterion} gives
\[
0\leq b_n\leq b_{n+1}\leq\E(\Lambda),
\qquad \lim_{n\to\infty}b_n=\E(\Lambda).
\]
Indeed, for every fixed $j$ and $n\geq j$,
$b_n\geq\ell_{j,n}\geq T_j-2^{-n}$. Therefore
$\lim_n b_n\geq\sup_jT_j=\E(\Lambda)$, while $b_n\leq\E(\Lambda)$.
This proves the assertion, including an infinite limit.
\end{proof}

\begin{remark}\label{rem:enumeration}
This construction evaluates individual tests and enumerates lower bounds;
it does not give a terminating decision procedure for arbitrary locally
finite sequences. The fixed-test error in \eqref{eq:finiteerror}
does not bound the tail of the supremum over all polynomials and radii.
The completeness proof supplies no a priori degree, support, or smoothing
bound for a requested approximation to $\E(\Lambda)$, and no convergence
rate for $b_n$ is asserted.

In particular, uniqueness of a sequence has the quantifier form
\[
\text{for every }M>0\text{ there exists a polynomial test with score exceeding }M.
\]
A finite computation can certify a particular lower bound. It cannot, by
itself, establish unboundedness, while finitely many small scores do not
certify that the supremum is finite. No general decidability or
undecidability result is claimed. The radial and lattice applications in
the preceding section are conclusive because they supply additional
analytic estimates.
\end{remark}

\medskip
\noindent\textbf{AI Disclosure.} Artificial-intelligence tools were used only
for language editing, LaTeX formatting, and limited assistance with local
mathematical reasoning. The proof strategy and mathematical development are
the authors' own work; all mathematical content was independently written and
checked by the authors, who take full responsibility for it.

\end{document}